\documentclass[10pt,twoside]{amsart}
\usepackage{amsmath}
\usepackage{amsfonts}
\usepackage{amssymb}
\usepackage{amsthm}

\theoremstyle{plain}
\newtheorem{thm}{Theorem}[section]
\newtheorem{cor}[thm]{Corollary}
\newtheorem{lem}[thm]{Lemma}
\newtheorem{prop}[thm]{Proposition}

\theoremstyle{remark}
\newtheorem{notation}[thm]{Notation}

\newcommand{\NZ}{\mbox{$\mathbb{N}$}}
\newcommand{\QZ}{\mbox{$\mathbb{Q}$}}
\newcommand{\CZ}{\mbox{$\mathbb{C}$}}
\newcommand{\Pone}{\mbox{$\mathbb{P}_1$}}
\newcommand{\Ptwo}{\mbox{$\mathbb{P}_2$}}
\newcommand{\Pthree}{\mbox{$\mathbb{P}_3$}}
\newcommand{\name}[1]{{\sc #1}\index{#1}}
\newcommand{\NEC}{\mbox{$\overline{NE(X)}$}\ }

\DeclareMathOperator{\Aut}{Aut}
\DeclareMathOperator{\Sing}{Sing}

\author{Stefan Kebekus}
\date{February 16, 1998}
\email{stefan.kebekus@uni-bayreuth.de}
\address{Stefan Kebekus\\ Mathematisches Institut der Universit\"at
Bayreuth\\ 95440 Bayreuth\\ Germany\\ FAX: +49 (0)921/55-2785 }
\title{Simple Models of Quasihomogeneous Projective 3-Folds}
\thanks{The author was supported by scholarships of the
Graduiertenkollegs ``Geometrie und mathematische Physik'' and
``Komplexe Mannigfaltigkeiten'' of the Deutsche
Forschungsgemeinschaft}

\begin{document} 

\begin{abstract}
Let $X$ be a projective complex 3-fold, quasihomogeneous with respect
to an action of a linear algebraic group. We show that $X$ is a
compactification of $SL_2/\Gamma$, $\Gamma$ a discrete subgroup, or
that $X$ can be equivariantly transformed into $\Pthree$, the quadric 
$\QZ_3$, or into certain quasihomogeneous bundles with very simple structure.

\end{abstract}

\maketitle
\tableofcontents

\section{Introduction}

Call a variety $X$ quasihomogeneous if there is a connected
algebraic group $G$ acting algebraically on $X$ with an open
orbit. Call a rational map $X\dasharrow Y$ equivariant if $G$ acts on
$Y$ and if the graph is stable under the induced action on $X\times
Y$. 

The class of varieties having an equivariant birational map to $X$ is
generally much smaller then the full birational equivalence class. The
minimal rational surfaces are good examples: they are all
quasihomogeneous with respect to an action of $SL_2$, but no two have
an $SL_2$-equivariant birational map between them. On the other hand,
if $X$ is any rational $SL_2$-surface, then the map to a minimal model
is always equivariant.

Generally, one may ask for a list of (minimal) varieties such that
every quasihomogeneous $X$ has an equivariant birational map to a
variety in this list.

We give an answer for $\dim X=3$ and $G$ linear algebraic:
\begin{thm}\label{mainthm}
Let $X$ be a 3-dimensional projective complex variety. Let $G$ be a
connected linear algebraic group acting algebraically and almost
transitively on $X$. Assume that the ineffectivity, i.e.~the kernel of
the map $G\rightarrow \Aut(X)$, is finite. Then either $G\cong SL_2$,
and $X$ is a compactification of $SL_2 / \Gamma$, where $\Gamma$ is
finite and not cyclic, or there exists an equivariant birational map
$X\dasharrow^{eq} Z$, where $Z$ is one of the following:
\begin{itemize}
\item $\Pthree$ or $\QZ_3$, the 3-dimensional quadric

\item a $\Ptwo$-bundle over $\Pone$ of the form
$\mathbb{P}(\mathcal{O}(e)\oplus\mathcal{O}(e)\oplus\mathcal{O})$.

\item a linear $\Pone$-bundle over a smooth quasihomogeneous surface
$Y$, i.e. $Z\cong \mathbb{P}(E)$, where $E$ is a rank-2 vector bundle
over $Y$. If $G$ is solvable, then $E$ can be chosen to be
split.
\end{itemize}
If $G$ is not solvable, then the map $X\dasharrow^{eq} Z$ factors into
a sequence $X\leftarrow \tilde{X}\rightarrow Z$, where the arrows
denote sequences of equivariant blow ups with smooth center.  
\end{thm}

A fine classification of the (relatively) minimal varieties involving
$SL_2$ will be given in a forthcoming paper.

The result presented here is contained the author's thesis. The author
would like to thank his advisor, Prof.~Huckleberry, and
Prof.~Peternell for support and valuable discussions.

\section{Existence of Extremal Contractions}\label{2sect}

The main tool we will use is \name{Mori}-theory. In order to utilize
it, we show that in our context extremal contraction always exist.

\begin{lem}
\label{exexcont}
Let $X$ and $G$ be as in~\ref{mainthm}, but allow for $\QZ$-factorial
terminal singularities. Then there exists a \name{Mori}-contraction. 
\end{lem}
\begin{proof}
Let $\pi:\tilde{X}\rightarrow X$ be an equivariant resolution of the
singularities of $X$. By $\tilde X$ being quasihomogeneous, we can
always find $\dim X$ elements $v_1, \ldots, v_n$ of the
\name{Lie}-algebra $Lie(G)$ such that the associated vector fields
$$
\tilde{v_i}(x) = \left. \frac{d}{dt}\right|_{t=0} \exp(tv_i) x
\quad\in H^0(\tilde{X},T\tilde{X})
$$
are linearly independent at generic points of $\tilde{X}$. In other
words, 
$$
\sigma := \tilde{v_1}\wedge \ldots \wedge \tilde{v_n} 
$$ is a non-trivial holomorphic section of the anticanonical bundle
$-K_{\tilde{X}}$. In effect, we have shown that $-K_{\tilde{X}}$ is
effective. Let $r$ be the index of $X$. By definition of terminal
singularities, the line bundle $rK_X$ is effective. Thus, we are
finished if we exclude the case that $rK_X$ is trivial.  Assume that
this is the case. The section $\sigma$ not vanishing on the smooth
points of $X$ implies that $X\setminus \Sing(X)$ is
$G$-homogeneous. But the terminal singularities are isolated. Thus, by
\cite[thm. 1 on p. 113]{H-Oe}, $X$ is a cone over a rational
homogeneous surface, a contradiction to $rK_X$
trivial.

In consequence $rK_X$ is effective and not trivial. So there is always
a curve $C$ intersecting an element of $|-K_X|$ transversally. Hence
$C.K_X <0$ and there must be an extremal contraction.
\end{proof}

\begin{cor}\label{correlm}
Let $X$ and $G$ be as in theorem~\ref{mainthm}, but allow for
$\QZ$-factorial terminal singularities. Let $\phi:X\rightarrow
Y$ an equivariant morphism with $\dim Y<3$. Then there is a relative
contraction over $Y$. 
\end{cor}
\begin{proof}
If $Y$ is a point, we are finished by lemma~\ref{exexcont}. Otherwise,
if $\eta \in Y$ generic, we know that the fiber $X_\eta$ is smooth,
does not intersect the singular set and is quasihomogeneous with
respect to the isotropy group $G_\eta$. So there exists a curve $C
\subset X_\eta$ with $C.K_{X_\eta} < 0$. Note that the adjunction
formula holds, since $X$ has isolated singularities and $X_\eta$ does
not intersect the singular set. Hence $K_{X_\eta} =
K_X|_{X_\eta}$, and there must be an extremal ray $C\subset \NEC$ such
that $\phi_*(C)=0$. Thus, there exists a relative contraction.
\end{proof}

Recall that all the steps of the \name{Mori} minimal model
program (i.e.~extremal contractions and flips) can be performed in an
equivariant way. For details, see \cite[chap. 3]{Ke96}.

\section{Equivariant Rational Fibrations}\label{pointchapter}

In this section we employ group-theoretical considerations in order to 
find equivariant rational maps from $X$ to varieties of lower
dimension. These will later be used to direct the minimal model program.

We start with the case that $G$ is solvable.

\begin{lem}\label{l11}
Let $X$ and $G$ be as in~\ref{mainthm}. Assume additionally that $G$
is solvable. Then there exists an equivariant rational map
$X\dasharrow^{eq}Y$ to a projective surface $Y$. 
\end{lem}
\begin{proof}
By $G$ being solvable, there exists a one-dimensional normal subgroup
$N$. Let $H$ be the isotropy group of a generic point, so that $\Omega
\cong G/H$, and consider the map 
$$
\Omega \cong G/H \rightarrow G/(N.H)
$$
Sine $N$ is not contained in $H$ (or else $G$ acted with positive
dimensional ineffectivity), the map has one-dimensional fibers. Now
$\dim G/(N.H) > 0$ and $G/(N.H)$ can always be equivariantly
compactified to a projective variety $Y$. This yields an equivariant
rational map $X\dasharrow^{eq} Y$. 
\end{proof}

Now consider the cases where $G$ is not solvable.

\begin{lem}\label{c0:4_1}
Let $X$ and $G$ be as above. Assume that $G$ is neither reductive nor
solvable. Then there exists an equivariant rational map
$X\dasharrow^{eq} Y$ such that either 
\begin{enumerate}
\item $Y \cong \Pthree$, and $X\dasharrow^{eq} Y$ is birational, or
$\dim Y =2$, or 
\item $\dim Y=1$, and there exists a normal unipotent group $A$ and a
semisimple group $S<G$, acting trivially on $Y$. The unipotent part
$A$ acts almost transitively on generic fibers.
\end{enumerate}
\end{lem}
\begin{proof}
Let $G=U\rtimes L$ be the  \name{Levi} decomposition of $G$, i.e. $U$ is
unipotent and $L$ reductive and define $A$ to be the center of
$U$. Note that $A$ is non-trivial. Since $A$ is canonically defined,
it is normalized by $L$, hence it is normal in $G$. Let $H$ be the
isotropy group of a generic point, $\Omega$ the open $G$-orbit, so
that $\Omega \cong G/H$, and consider the map 
$$
\Omega \cong G/H \rightarrow G/(A.H)
$$
There are two things to note. The first is that $A$ is not contained in $H$
(or else $G$ acted with positive dimensional ineffectivity). 
So $\dim G/(A.H) < 3$. If $\dim G/(A.H) > 0$, it can always be equivariantly 
compactified $G/(A.H)$ to a variety $Y$ yielding an equivariant rational
map $X\dasharrow Y$. If $\dim G/(A.H)=2$, we can stop here. If $\dim
G/(A.H)=1$, then note that $A$ acts transitively on the fiber
$A.H/H$. If $A.H$ does not contain a semi-simple group, we argue as in
lemma~\ref{l11} to find a subgroup $H'$, $H<H'<A.H$ such that $\dim
H'/H =1$. Then $\dim G/H'=2$, and again we are finished.

If $\dim G/(A.H) = 0$, then $A$ acts transitively on $\Omega$.  In
this case $A\cong \CZ^n$, and hence (because the $G$-action is
algebraic) $\Omega \cong \CZ^3$. The theorem on \name{Mostow}
fibration (see e.g. \cite[p. 641]{Hei91}) yields that $L$ has to have
a fixed point in $\Omega$. Therefore, without loss of generality,
$L<H$. As a next step, consider the group $B:=(U\cap H)^0$. Since both
$U$ and $H$ are normalized by $L$, $B$ is as well. Elements in $A$
commute with all elements of $U$, hence $A.B$ normalizes $B$ as
well. Then $B$ is a normal subgroup of $U\rtimes L=G$. Note that $A.B
=U$, because $A.B = A.(H\cap U) = (A.H)\cap U = G \cap U =
U$. Consequently $B$ does acts trivially and so it is trivial.

We are now in a position where we may write $G=A\rtimes_{\rho} L$,
where $\rho$ is the action of $L$ on $A$ ($L$ acting by
conjugation). Now $H=L$, hence $A\cong \Omega \cong \CZ^3$ and the
$L$-action on $A\cong (\CZ^3,+)$ has to be linear. So $G$ is a
subgroup of the affine group and $\Omega$ can be equivariantly
compactified to $\Pthree$, yielding an equivariant rational map
$X\dasharrow^{eq} \Pthree$.
\end{proof}

We study case (1) of the preceding proposition in more detail.

\begin{lem}
Let $X$ be as above and assume that $G$ be reductive. Assume
furthermore that $G$ is not semisimple. Then there is an equivariant
rational map $X \dasharrow^{eq} Z$, where $\dim Z =2$.
\end{lem}

\begin{proof}
As a first step, recall that $G = T.S$, where $S$ semisimple, $T$ a
torus, and $S$ and $T$ commute and have only finite intersection. 
If $\eta$ is a point in the open orbit and $G_{\eta}$ the
associated isotropy group, then $T \not \subset G_{\eta}$, or
otherwise  $T$ would not act at all. For that
reason we will be able to find a 1-parameter group $T_1 < T$, $T_1
\not \subset G_{\eta}$ and consider the map 
$$
\Omega := G/G_{\eta} \rightarrow G/(T_1.G_{\eta}).
$$
Since $T_1$ has non-trivial orbits, $\dim G/(T_1.G_{\eta}) =2$. If we 
compactify the latter in an equivariant way to a variety
$Z$, we automatically obtain a an equivariant rational map $X
\dasharrow^{eq} Z$ as claimed.
\end{proof}

\begin{lem}\label{maps2p2}
Suppose $G$ is semisimple. Then one of the following holds:
\begin{enumerate}
\item $G\cong SL_2$ and the open orbit $\Omega$ is isomorphic to
$SL_2/\Gamma$, where $\Gamma$ is discrete and not contained in a
\name{Borel} subgroup.

\item $X \cong \Pthree$

\item $X$ is isomorphic to $F_{1,2}(3)$, the full flag variety

\item $X$ is homogeneous and either $X\cong \QZ_3$, the 3-dimensional
quadric or $X$ is a direct product involving only $\Pone$ and $\Ptwo$.

\item $X$ admits an equivariant rational map $X\dasharrow^{eq} Y$ onto a
surface.
\end{enumerate}
\end{lem}
\begin{proof}
If $G\cong SL_2$, and $\Gamma$ is embeddable into a \name{Borel} group
$B$, then $\Gamma$ is in fact embeddable into a 1-dimensional torus
$T$. Consider the map $G/\Gamma \rightarrow G/T$, and 
we are finished. 

Assume for the rest of this proof that $G \not
\cong SL_2$. Then the claim is already true in the complex analytic
category: see \cite[p. 3]{Winkelmann}. One must exclude torus bundles
by the fact that they never allow an algebraic action of a {\em linear}
algebraic group.
\end{proof}

We summarize a partial result:

\begin{cor}\label{mappingcor}
Let $X$ and $G$ be as above. If there exists an equivariant map
$X\dasharrow^{eq} \Pone$ and no such map to $\Pthree$ or to a surface,
then $G$ is not solvable and there exist subgroups $S$ and $A$ as in
lemma~\ref{c0:4_1}.
\end{cor}

\section{The case that $Y$ is a curve}
\label{curvechap}

In this section we investigate relatively minimal models over
$\Pone$. The main proposition is:

\begin{prop}
\label{MP2}Let $X$ and $G$ be as in~\ref{mainthm}, but allow for
$\QZ$-factorial terminal singularities. Assume that $\phi:X\rightarrow
\Pone$ is an extremal contraction. Assume additionally that there does
not exist an equivariant rational map $X\dasharrow^{eq} Y$, where
$\dim Y = 2$ or $Y\cong \Pthree$. Then
$$
X \cong
\mathbb{P}(\mathcal{O}_{\Pone}(e) \oplus \mathcal{O}_{\Pone}(e) \oplus
\mathcal{O}_{\Pone}),
$$
with $e>0$. In particular, $X$ is smooth.
\end{prop}

\begin{proof}
As a first step, we show that the generic fiber $X_\eta$ is isomorphic
to $\Ptwo$. As $\phi$ is a \name{Mori}-contraction, $X_\eta$ is a
smooth \name{Fano} surface. By corollary~\ref{mappingcor}, the
stabilizer $G_\eta<G$ of $X_\eta$ contains a unipotent group $A$
acting almost transitively on $X_\eta$ and a semisimple part $S$. This
already rules out all \name{Fano} surfaces other than
$\Ptwo$. Furthermore, $S\cong SL_2$.  Note that $G_\eta$ stabilizes a
unique line $L\subset X_\eta$ and that $S$ acts transitively on $L$.

Set $D' := \overline{G.L}$ and remark that $D'$ intersects the generic
$\phi$-fiber in the unique $G_\eta$-stable line: $D'\cap X_\eta =
L$. We claim that $D'$ is \name{Cartier}. The desingularization
$\tilde{D}'$ has a map to $\Pone$, the generic fiber is isomorphic to
$\Pone$ and $S$ acts on all the fibers. Thus, $\tilde{D}'$ is
isomorphic to $\Pone\times\Pone$, and $S$ does not have a fixed point
on $D'$. Consequently, $\tilde{D}'$ does not intersect the singular
set of $X$ and is \name{Cartier}.

Take $D''$ to be an ample divisor on $Y$. As $\phi$ is a
\name{Mori}-contraction, the line bundle $L$ associated to
$D:=D'+n\phi^*(D'')$, $n>>0$, is ample on $X$. In this setting, a
theorem of \name{Fujita} (cf. \cite[Prop. 3.2.1]{Belt95}) yields that
$X$ is of the form $\mathbb P(E)$, where $E$ is a vector bundle on $\Pone$. 

The transition functions of $E$ must commute with $S$, but the only
matrices commuting with $SL_2$ are $Diag(\lambda,\lambda,\mu)$, hence
$E=\mathcal{O}(e)\oplus\mathcal{O}(e)\oplus\mathcal{O}(f)$ and $X
\cong
\mathbb{P}(\mathcal{O}(e-f)\oplus\mathcal{O}(e-f)\oplus\mathcal{O})$. 
\end{proof}

For future use, we note

\begin{lem}\label{lnof}
Let $X$ and $G$ be as in proposition~\ref{MP2}. Then, by
equivariantly blowing up and down, $X \dasharrow^{eq}
\mathbb{P}(\mathcal{O}(e')\oplus\mathcal{O}(e')\oplus\mathcal{O})$
where the latter does not contain a $G$-fixed point.
\end{lem}
\begin{proof}
The semisimple group $S$ fixes a unique point of each $\phi$-fiber, so
that there exists a curve $C$ of $S$-fixed points. Suppose that $G$ has a
fixed point $f$. Then $f\in C$, and we can perform an elementary
transformation $X\dasharrow^{eq} X'$ with center $f$, i.e. if $X_\mu$
is the $\phi$-fiber containing $f$, then we blow up $f$ and blow down
the strict transform of the $X_\mu$, again obtaining a linear
$\Ptwo$-bundle of type $\mathbb{P}(\mathcal{O}(e) \oplus
\mathcal{O}(e) \oplus \mathcal{O})$. This transformation exists, as
has been shown in \cite{Maru73}. Since all the centers of the blow-up
and -down are $G$-stable, the transformation is equivariant.

We will use this transformation in order to remove $G$-fixed points. Let
$g\in G$ be an element not stabilizing $C$. The curves $gC$ and $C$
meet in $f$. We know that after finitely many blow-ups of the
intersection points of $C$ and $gC$, the curves become disjoint, so
that there no longer exists a $G$-fixed point! This, however, is
exactly what we do when applying the elementary transformation.
\end{proof}

\section{The case that $Y$ is a surface}\label{surfchap}

The cases that $G$ is solvable or not solvable behave differently in
many respects. Here we have to treat them separately.

\subsection{The case $G$ solvable}

We will show that in this situation the open $G$-orbit can be
compactified in a particularly simple way.

\begin{prop}
\label{EQMap_to_splitting} Let $X$ and $G$ be as in
theorem~\ref{mainthm}. Assume additionally 
that $G$ is solvable and $\phi:X\rightarrow Y$ is an equivariant map
with connected fibers onto a smooth surface. Then there exists a
splitting rank-2 vector bundle $E$ on $Y$ and an equivariant
birational map $X\dasharrow^{eq} \mathbb{P}(E)$.
\end{prop}

We remark that if $y\in Y$ is contained in the open $G$-orbit, then
it's preimage is quasihomogeneous with respect to the isotropy group
$G_y$, hence isomorphic to $\Pone$. As a first step in the proof of
proposition~\ref{EQMap_to_splitting}, we show the existence of very
special divisors in $X$, namely:

\begin{notation}
We call a divisor $D\subset X$ a ``rational section'' if it intersects
the generic $\phi$-fiber with multiplicity one.
\end{notation}

In our context, such divisors always exist:

\begin{lem}\label{l23}
Let $\phi:X\rightarrow Y$ be as in lemma~\ref{EQMap_to_splitting} and assume
additionally that there exists a group $H^*\cong \CZ^*$, acting
trivially on $Y$. Let $D_X'$ be the fixed point set of the
$H^*$-action. Then $D_X'$ contains two rational sections as
irreducible components.
\end{lem}
\begin{proof}
Let $D_X$ be the union of those irreducible divisors in $D'_X$ which
are not preimages of curves or points by $\phi$. The subvariety $D_X$
intersects every generic $\phi$-fiber at least once. Hence $D_X \not =
0$.

We claim that the set of branch points $$ M:= \{ y\in Y:
\#(\phi^{-1}(y) \cap D_X) = 1 \} $$ is discrete. Linearization of the
$H^*$-action yields that for any point $f \in D_X \setminus \Sing(X)$,
there is a unique $H^*$-stable curve intersecting $D_X$ at
$f$. Furthermore, the intersection is transversal. Assume $\dim M\geq
1$ and let $y$ be a generic point in $M$. Then $\dim \phi^{-1}(y)=1$
and $\phi^{-1}(y)=1$ contains a smooth curve $C$ as an irreducible
component intersecting $D_X$. Now $C.D_X =1$, and because $C\cap D_X$
was the only intersection point by assumption, $\phi^{-1}(y).D_X=1$,
contradicting $D_X$ intersecting the generic $\phi$-fiber twice.

Set 
$$
N:=\{\mu \in Y | \dim (X_\mu \cap D_X) > 0 \}\cup M\cup
\phi(\Sing(X)).
$$ By definition, $N$ is finite and $D_X$ is a 2-sheeted cover over $Y
\setminus N$. Now $Y$ is smooth and quasihomogeneous with respect to
an algebraic action of the linear algebraic group $G$, hence rational.
This implies that $Y \setminus N$ is simply connected. Hence $D_X$ has
two connected components over $Y \setminus N$. Now the set
$D_X\cap\phi^{-1}(N)$ is just a curve. Therefore $D_X$ cannot be
irreducible.
\end{proof}

\begin{lem}\label{224}
Under the assumptions of lemma~\ref{EQMap_to_splitting}, 
there exists a $G$-stable rational section $E_1 \subset X$. 
\end{lem}
\begin{proof}
If $G$ is a torus, then there exists a subgroup $T_1$ acting trivially
on $Y$. In this case we are finished by applying
lemma~\ref{l23}. Hence, assume that the unipotent part $U$ of $G$ is
non-trivial. Let $\eta \in Y$ be a generic point and $x\in X_\eta
\setminus \Omega$, where $\Omega$ denotes the open $G$-orbit in
$X$. If $x$ is unique, then we are done by setting
$E_1:=\overline{G.x}$. Similarly, if $U$ acts almost transitively on
$Y$, then it's isotropy at $\eta$ is connected and we may set
$E_1:=\overline{U.x}$.

If neither holds, then necessarily $\dim U=1$, and we can assume that
$U$ acts non-trivially on $Y$, or else $X_\eta \setminus \Omega$
consists of a single point and we are finished as above. Let $T_1$ be
a 1-dimensional subgroup of a maximal torus such that $I:=U.T_1$ acts
almost transitively on $Y$. If $\eta \in Y$ is generic, the isotropy
group $I_\eta$ is cyclic: $I_\eta$ has two fixed points in
$X_\eta$. Consequently, there exist at least two $I$-orbits whose
closures $D_i$ are rational sections. 

Note that $I$ is normal in $G$, i.e.~all elements of $G$ map
$I$-orbits to $I$-orbits. If $D_i$ are the only rational sections
occurring as closures of $I$-orbits, they are automatically
$G$-stable. Otherwise, all $I$-orbits are mapped injectively to $Y$,
and at least one of these is $G$-stable.
\end{proof}

The existence of $E_1$ already yields a map to a $\Pone$-bundle.

\begin{lem}\label{230}
Under the assumptions of lemma~\ref{EQMap_to_splitting}, 
there exists a rank-2 vector bundle $E$ on $Y$ (not necessarily split)
and an equivariant birational map $X\dasharrow^{eq} \mathbb{P}(E)$.
\end{lem}
\begin{proof}
Set $E := (\phi_*(\mathcal{O}_X(E)))^{**}$. As a reflexive sheaf is
locally free on a smooth surface, $E$ is a vector bundle. If $\Omega_Y
\subset Y$ is the open orbit, $\phi^{-1}(\Omega_Y)\cong
\mathbb{P}(E|_{\Omega_Y})$ (cf. \cite[Prop. 3.2.1]{Belt95}), inducing
a birational map $\psi:X\dasharrow \mathbb{P}(E)$. Note that
$\phi_*(\mathcal{O}_X(E))$ is torsion free. In particular,
$\phi_*(\mathcal{O}_X(E))$ is locally free over a $G$-stable cofinite
set $Y_0 \subset Y$ so that, by the universal property of $Proj$,
$\psi$ is regular over $Y_0$. As $\psi|_{Y_0}$ is proper, it is
equivariant. The automorphisms over $Y_0$ extend to the whole of
$\mathbb{P}(E)$ by the \name{Riemann} extension theorem. Hence $\psi$
is equivariant as claimed. 
\end{proof}

In order to show that $E$ can be chosen to be split, we need to find
another rational section. We will frequently deal with the following
situation, for which we fix some notation.

\begin{notation}\label{Fnot}
Let $\phi:X\rightarrow Y$ be as above and assume that
there exists a map $\pi: Y\rightarrow Z\cong\Pone$, e.g. if $Y$ is
isomorphic to a (blown-up) \name{Hirzebruch} surface $\Sigma_n$. Then,
if $F \in Z$ is a generic point, set $F_Y:=\pi^{-1}(F)$ and
$F_X:=\phi^{-1}(F_Y)$. 
\end{notation}

\begin{lem}
In the setting of proposition~\ref{EQMap_to_splitting}, there exists a
second rational section $E_2$. If $E_1$ is as constructed in
lemma~\ref{224}, then $E_1 \cap E_2$ is $G$-stable.
\end{lem}
\begin{proof}
If $G$ is a torus, we are finished, as we have seen in the proof of
lemma~\ref{224}. Hence we may assume that $\dim U>0$, where $U$ is the
unipotent part of $G$.  

Suppose that $U$ acts trivially on $Y$. Then we are able to choose a
2-dimensional torus $T<G$ such that $T$ acts almost transitively on
$Y$. If $\eta\in Y$ is generic, then the isotropy group $T_\eta$ may
not be cyclic, but since it has to fix the unique $U$-fixed point in
$X_\eta$, its image $T_\eta\rightarrow \Aut(X_\eta)$ is contained in
a \name{Borel} group, hence cyclic. Consequently, $T_\eta$ fixes
another point $x$, and we may set $E_2:=\overline{T.x}$. 

The other case is that $U$ acts non-trivially on $Y$. We need to
consider a mapping $\pi:Y\rightarrow Z\cong \Pone$. If $Y\cong \Sigma_n$, or
a blow-up, there is no problem. If $Y\cong \Ptwo$, we note that,
by $G$ being solvable and \name{Borel}'s fixed point theorem (see
\cite[p. 32]{H-Oe}), there exists a $G$-fixed point $y \in Y$. We can
always blow up $y$ and $X_y$ in order to obtain a new $\Pone$-bundle
over $\Sigma_1$. If we are able to construct our rational sections
here, then we can simply take their images to be the desired rational
sections in the variety we started with. So let us assume that $Y \not
\cong \Ptwo$. 

There exists a 1-dimensional normal unipotent subgroup $U_1<G$. Assume
first that $U_1$ acts non-trivially on $Z$. Using
notation~\ref{Fnot}, $F_Y$ is isomorphic to $\Pone$, $F_X$ to a 
\name{Hirzebruch} surface $\Sigma_n$. Choose a section $\sigma \subset
F_X$ with the property that $\phi(\sigma \cap E_1)$ does not meet the
open $G$-orbit in $Y$. As the stabilizer of $F_X$ in $G$ stabilizes
$E_1$, so that $E_1\cap F_X$ is either the infinity- or zero-section
in $F_X \cong \Sigma_n$ or the diagonal in $F_X \cong \Sigma_0$, and
$G$ stabilizes a section of $Y\rightarrow \Pone$, this can always be
accomplished. Set $E_1:=\overline{U_1.\sigma}$.

Secondly, we must consider the case that $U_1$ acts trivially on
$Z$. We proceed similarly to the above. Choose a 1-dimensional group
$G_1 <G$ such that the $G_1$-orbit in $Z$ coincides with that of
$G$. Now $G_1$ stabilizes at least one section $\sigma_Y \subset Y$ over
$Z$ which is not $U_1$-stable! Set $\sigma_X := \phi^{-1}(\sigma_Y)$
and consider a section $\sigma \subset \sigma_X$ over $\sigma_Y$ 
such that $\phi(\sigma \cap E_1)$ is disjoint from the open $G$-orbit
in $Y$. Then $E_1 := \overline{U_1.\sigma}$ is the divisor we were
looking for. 

We still have to show that the intersection $E_1 \cap E_2$ is
$G$-stable. Note that by construction, $\phi(E_1\cap E_2)$ does not
meet the open $G$-orbit in $Y$. This, together with $E_1$ being
$G$-stable, yields the claim. 
\end{proof}

We shall use the second rational section in order to transform $E$
into a splitting bundle.

\subsubsection{Eliminating vertical curves}

If $S \subset \phi(E_1\cap E_2)$ is an irreducible curve which is a
$\phi$-fiber, say that $E_1$ and $E_2$ intersect vertically in $S$. We
know that after blowing up $S$ we obtain a $\Pone$-bundle over the
blow-up of $Y$. Furthermore, the process is equivariant. The proper
transforms of $E_1$ and $E_2$ are again rational sections. If they
still intersect vertically, the blow-up procedure can be applied
again. So we eventually get a sequence of blow-ups. The strict
transforms of the $E_1$ and $E_2$ are again rational sections in
$X_i$. We denote them by $E^{i}_1$ or $E^{i}_2$, respectively. By the
theorem on embedded resolution, we have:

\begin{lem}
\label{df3}
The sequence described above terminates, i.e. there
exists a number $i\in \NZ$ such that the strict transforms $E^i_1$ and
$E^i_2$ do not intersect vertically.
\end{lem}

\subsubsection{Eliminating horizontal curves}

We may now assume that $E_1$ and $E_2$ do not intersect vertically. 
Let $S \subset \phi(E_1\cap E_2)$ be an irreducible curve. Then $S$
gives rise to an elementary transformation as ensured by
\cite{Maru73}. Again, the transformation is equivariant and the strict transforms of $E_1$ and $E_2$ are
rational sections. If they still intersect over $S$, we transform as
before. Again one may use the embedded resolution to show
(cf. \cite[thm. 5.30]{Ke96} for details): 

\begin{lem}
\label{T_helps}
The sequence described above terminates after finitely many
transformations, i.e. there exists a $j \in \NZ$ such that for all
curves $C \in E^{(j)}_1 \cap E^{(j)}_2$ it follows that $\phi^{(j)}(C)
\not = S$. Furthermore, if $E_1$ and $E_2$ do not intersect
vertically, then $E^{(i)}_1$ and $E^{(i)}_2$ do not intersect 
vertically for all $i$.
\end{lem}

\subsubsection{The construction of independent sections}

By lemma~\ref{df3} the variety $X$ can be transformed into
a $\Pone$-bundle such that the strict transforms of $E_1$ and $E_2$
do not intersect in fibers. A second transformation will rid us of
curves in $E_1 \cap E_2$ which are not contained in fibers. Since the
latter transformation  does not create new curves in the
intersection, the strict transforms of $E_1$ and $E_2$ eventually
become disjoint. The resulting space is the compactification of a line
bundle.

\begin{lem}
If $E_1$ and $E_2$ do not intersect, $X$ is the compactification of a
line bundle. 
\end{lem}
\begin{proof}
Since $E_1$ and $E_2$ are disjoint, neither contains a fiber. Thus they
are sections.
\end{proof}

As a net result, we have shown proposition \ref{EQMap_to_splitting}.

\subsection{The case $G$ not solvable}

As first step, we show that $X$ is again a linear $\Pone$-bundle. We
do this under an additional hypothesis which will not impose problems
in the course of the proof of theorem~\ref{mainthm}.

\begin{lem}\label{231}
Let $X$ and $G$ be as in theorem~\ref{mainthm}, but allow for
$\QZ$-factorial terminal singularities. Let $\phi:X\rightarrow Y$
a \name{Mori}-contraction to a surface and assume additionally
that $G$ is not solvable and that there exists an equivariant morphism
$\psi:Y\rightarrow Y'$, where $Y'$ is a smooth surface. Then $X$ and
$Y$ are smooth and $X$ is a linear $\Pone$-bundle over $Y$.
\end{lem}
\begin{proof}
First, show that all $\phi$-fibers are of dimension 1. If there exists
a fiber $X_\mu$ which is not 1-dimensional, then $\dim X_\mu=2$. Take
a curve $C\subset Y$ so that $\mu \in C$. Set $D:=
\overline{\phi^{-1}(C\setminus \mu)}$. The divisor $D$ intersects an
irreducible component of $X_\mu$. Now take a curve $R \subset X_\mu$
intersecting $D$ in finitely many points. We have $R.D > 0$. However,
all generic $q$-fibers $X_\eta$ are homologeously equivalent to $R$
(up to positive multiples). So $X_\eta.D > 0$, contradicting the
definition of $D$.

Secondly, we claim that $X$ is smooth. Assume to the contrary and let
$x\in X$ be a singular point, $\mu:= \phi(x)$. Recall that terminal
singularities in 3-dimensional varieties are isolated. Thus, if $S$ is
the semisimple part of $G$, then the fiber $X_\mu$ through $x$ is
pointwise $S$-fixed. Linearizing the $S$-action at a generic point
$y\in X_\mu$, the complete reducibility of the $S$-representation
yields an $S$-quasihomogeneous divisor $D$ which intersects $X_\mu$
transversally in $y$ and is \name{Cartier} in a neighborhood of
$y$. The induced map $D\rightarrow Y'$ must be unbranched: $Y'$
contains an $S$-fixed point and is therefore isomorphic to $\Ptwo$;
but there is no equivariant cover of this other than the identity. So
$D$ is a rational section which is \name{Cartier} over a neighborhood
of $\mu$. If $H\in Pic(Y)$ is sufficiently ample, then $D+\phi^*(H)$
is ample, and \cite[Prop. 3.2.1]{Belt95} applies, contradicting $X$
singular.

As $X$ is smooth, the same theorem shows that in order to prove the
lemma it is sufficient to show that there exists a rational
section. If all the simple factors of $S$ have orbits of dimension
$\leq 2$, then, after replacing the factors by their \name{Borel}
groups, we obtain a solvable group $G'$, acting almost transitively as
well. In this case lemma~\ref{224} applies.

If $S'<S$ is a simple factor acting with 3-dimensional orbit on $X$,
its action on $Y$ is almost transitively. In particular, there exists
a 2-dimensional group $B<S$, isomorphic to a \name{Borel} group in
$SL_2$, acting almost transitively on $Y$, too. As in the proof of
lemma~\ref{224}, $B$ has cyclic isotropy at a generic point of $Y$ and
so there exist two rational sections which are compactifications of
$B$-orbits. 
\end{proof}

\section{Proof of theorem~\ref{mainthm}}\label{link}

Prior to proving theorem~\ref{mainthm}, we still need to describe
equivariant maps to $\Pthree$ in more detail:

\begin{lem}
\label{P3cases}
Let $X\dasharrow^{eq} \Pthree$ be an equivariant birational map. Then
either $X$ has an equivariant rational fibration with 2-dimensional
base variety or $X$ and $\Pthree$ are equivariantly linked by a
sequence of blowing ups of $X$ followed by a sequence of blow-downs.
\end{lem}
\begin{proof}
If the $G$-action on $\Pthree$ has a fixed point, we can blow up this
point and obtain a map from the blown-up $\Pthree$ to $\Ptwo$. If
there is no such $G$-fixed point in $\Pthree$, then after replacing
$X$ by an equivariant blow-up, there is a regular equivariant map
$\phi: X\rightarrow \Pthree$. Recall that such a map factors through an
extremal contraction. As the base does not have a fixed point, the
classification of extremal contractions of smooth varieties yields the
claim.
\end{proof}

Now we compiled all the results needed to finish the

\begin{proof}[Proof of theorem~\ref{mainthm}]
Given $X$, we apply lemmata~\ref{l11}--\ref{maps2p2}. Unless $X\cong
Q_3$, $F_{1,2}(3)$ or a compactification of $SL_2/\Gamma$, $\Gamma$
not cyclic, there exists an equivariant map $X\dasharrow^{eq} Y$,
where $Y$ is smooth and $Y\cong \Pthree$, $\dim (Y)=2$ or, if no other
case applies, $\dim (Y)=1$.

If $Y\cong \Pthree$, then, by lemma~\ref{P3cases}, we may replace
$\Pthree$ by a surface, or else we are finished. 

In the case of a map to $Y$ with $\dim Y<3$, we can blow up $X$
equivariantly so as to have a morphism $\tilde{X}\rightarrow
Y$. Recalling that all steps in the minimal model program
(i.e.~contractions and flips) are equivariant, we may perform a
relative minimal model program over $Y$. In this situation
corollary~\ref{correlm} yields that the program does not stop unless
we encounter a contraction of fiber type $X'\rightarrow Y'$ and $\dim
Y'<3$. Note that $\dim Y' \geq \dim Y$.

In case that $Y'$ is a surface, $X'$ is the projectivization of a line
bundle or can be equivariantly transformed into one
(cf. lemma~\ref{230} and \ref{231}). If $G$ is solvable,
proposition~\ref{EQMap_to_splitting} allows us to transform $X$ into
the projectivization of a splitting bundle over a surface.

If $\dim Y'=1$, and there does not exist a map to one of the other
cases, $X\cong \mathbb{P}(\mathcal{O}(e) \oplus \mathcal{O}(e) \oplus
\mathcal{O})$ over $\Pone$, as was shown in proposition~\ref{MP2}. 

We still have to show that if $G$ is not solvable, the map to one of
the models in our list factors into equivariant monoidal
transformations. Recall that it suffices to show that, after equivariantly
blowing up, if necessary, the minimal models do not have a $G$-fixed
point. We do a case-by-case checking:
\begin{description}
\item[$\Ptwo$-bundles over $\Pone$] By lemma~\ref{lnof}, these can be
chosen not to contain a fixed point.

\item[$\Pone$-bundles over a surface $Y$] If $S$, the semisimple part
of $G$, acts trivially on $Y$, we can stop. Otherwise, if the
$S$-action on $Y$ has a fixed point $f$, we blow up $f$ and the fiber
over $f$, obtaining a $\Pone$-bundle over $\Sigma_1$. Recall that
actions of semisimple groups on $\Sigma_n$ never have fixed points.

\item[$\Pthree$] This case has already been handled in
lemma~\ref{P3cases}.

\item[$\boldsymbol{SL_2/\Gamma}$] after desingularizing and blowing up
all fixed points, if any, the compactification of $SL_2/\Gamma$ is
fixed point free. Otherwise, linearization at a fixed point yields a
contradiction to $S$ acting almost transitively.

\item[other cases] The other cases occur only when homogeneous 
(cf.~lemma~\ref{maps2p2}).
\end{description}
This shows the claim.
\end{proof}

\end{document}